\documentclass[11pt,a4paper]{article}
\usepackage[T1]{fontenc}
\usepackage[utf8]{inputenc}
\usepackage{authblk}
\usepackage{geometry}
\author{Yue Feng}
\affil{Department of Mathematics, National University of Singapore, Singapore 119076, Singapore}

\date{} 
\usepackage{amssymb,amsfonts,amsmath}
\usepackage{graphicx}
\usepackage{color}
\usepackage{booktabs}
\usepackage{multirow}
\usepackage{threeparttable}
\newtheorem{theorem}{Theorem}
\newtheorem{lemma}{Lemma}
\newtheorem{remark}{Remark}
\numberwithin{equation}{section}
\numberwithin{theorem}{section}
\numberwithin{lemma}{section}
\numberwithin{remark}{section}

\graphicspath{{figures/}}

\begin{document}

\title{Long time error analysis of the fourth-order compact finite difference methods for the nonlinear Klein-Gordon equation with weak nonlinearity}

\maketitle

\begin{abstract}
We present the fourth-order compact finite difference (4cFD) discretizations for the long time dynamics of the nonlinear Klein-Gordon equation (NKGE), while the nonlinearity strength is characterized by $\varepsilon^p$ with a constant $p \in \mathbb{N}^+$ and a dimensionless parameter $\varepsilon \in (0, 1]$. Based on analytical results of the life-span of the solution, rigorous error bounds of the 4cFD methods are carried out up to the time at $O(\varepsilon^{-p})$. We pay particular attention to how error bounds depend explicitly on the mesh size $h$ and time step $\tau$ as well as the small parameter $\varepsilon \in (0, 1]$, which indicate that, in order to obtain `correct' numerical solutions up to the time at $O(\varepsilon^{-p})$, the $\varepsilon$-scalability (or meshing strategy requirement) of the 4cFD methods should be taken as: $h = O(\varepsilon^{p/4})$ and $\tau = O(\varepsilon^{p/2})$. It has better spatial resolution capacity than the classical second order central difference methods. By a rescaling in time, it is equivalent to an oscillatory NKGE whose solution propagates waves with wavelength at $O(1)$ in space and $O(\varepsilon^p)$ in time. It is straightforward to get the error bounds of the oscillatory NKGE in the fixed time. Finally, numerical results are provided to confirm our theoretical analysis.
\end{abstract}


\section{Introduction}
\label{intro}
The nonlinear Klein-Gordon equation (NKGE) is a relativistic (and nonlinear) version of the Schr\"{o}dinger equation and widely used to describe the motion of a spinless particle \cite{BCZ,BD,DEGM,SJJ,M}. This equation has gained much attention in nonlinear optics, solid state physics and quantum field theory\cite{BS,DG,W}. We consider the following NKGE on a torus $\mathbb{T}^d (d=1,2,3)$ \cite{FV,SV} 

\begin{equation}
\begin{split}
&\partial_{tt} u({\bf{x}},t) -\Delta u({\bf{x}},t) +  u({\bf{x}},t) + \varepsilon^p u^{p+1}({\bf{x}},t)= 0, \quad  {\bf{x}} \in \mathbb{T}^d,\quad t > 0,\\
&u({\bf{x}}, 0) = \phi({\bf{x}}),\quad \partial_t u( {\bf{x}}, 0) = \gamma( {\bf{x}}),\quad  {\bf{x}} \in \mathbb{T}^d,
\end{split}
\label{eq:WNE}
\end{equation}
with $p \in \mathbb{N}^+$. Here $t$ is time, ${\bf{x}} \in \mathbb{T}^d$ is the spatial coordinates, $u: = u({\bf{x}},t)$ is a real-valued scalar field, $0 < \varepsilon \leq 1 $ is a dimensionless parameter, and $\phi( {\bf{x}})$ and $\gamma( {\bf{x}})$ are two given real-valued functions which are independent of $\varepsilon$. The NKGE \eqref{eq:WNE} is time symmetric or time reversible. In addition, if $u(\cdot, t) \in H^1(\mathbb{T}^d)$ and $\partial_t u(\cdot, t) \in L^2(\mathbb{T}^d)$, it also conserves the {\sl energy} \cite{BD,DXZ}, i.e.,
\begin{equation}
\begin{split}
E(t) := & \int_{\mathbb{T}^d} \left[ |\partial_t u ({\bf{x}}, t)|^2 + |\nabla u({\bf{x}}, t)|^2 + |u({\bf{x}}, t)|^2 +\frac{2\varepsilon^p}{p+2} |u({\bf{x}}, t)|^{p+2}  \right] d {\bf{x}} \\
 \equiv & \int_{\mathbb{T}^d} \left[ |\gamma({\bf{x}})|^2 + |\nabla \phi({\bf{x}})|^2 + |\phi({\bf{x}})|^2 +\frac{2\varepsilon^p}{p+2} |\phi({\bf{x}})|^{p+2}  \right] d {\bf{x}} \\
 = & E(0)=O(1), \quad t \geq 0.
\end{split}
\label{eq:Energy_u}
\end{equation}

When $0 < \varepsilon \ll1 $, introducing $w({\bf{x}},t)  = \varepsilon u({\bf{x}},t)$, the NKGE \eqref{eq:WNE} with weak nonlinearity can be reformulated as the following NKGE  with $O(\varepsilon)$ initial data:
\begin{equation}
\begin{split}
&\partial_{tt} w({\bf{x}},t)  -\Delta w({\bf{x}},t)  +   w({\bf{x}},t)  + w^{p+1}({\bf{x}},t) = 0, \quad {\bf{x}} \in \mathbb{T}^d,\quad t > 0,\\
&w({\bf{x}}, 0) = \varepsilon \phi( {\bf{x}}),\quad \partial_t w({\bf{x}}, 0) = \varepsilon \gamma({\bf{x}}),\quad {\bf{x}} \in \mathbb{T}^d.
\label{eq:SIE}
\end{split}
\end{equation}
Again, the above NKGE  (\ref{eq:SIE}) is time symmetric or time reversible and conserves the {\sl energy}, i.e.,
\begin{equation}
\begin{split}
\bar{E}(t) := &\int_{\mathbb{T}^d} \left[ |\partial_t w ({\bf{x}}, t)|^2 + |\nabla w({\bf{x}}, t)|^2 + |w({\bf{x}}, t)|^2 +\frac{2}{p+2} |w({\bf{x}}, t)|^{p+2}  \right] d {\bf{x}} =\varepsilon^2 E(t)\\
\equiv & \int_{\mathbb{T}^d} \left[ \varepsilon^2 |\gamma({\bf{x}})|^2 + \varepsilon^2 |\nabla \phi({\bf{x}})|^2 + \varepsilon^2|\phi({\bf{x}})|^2 +\frac{2\varepsilon^{p+2}}{p+2} |\phi({\bf{x}})|^{p+2}  \right] d {\bf{x}} \\
= & \bar{E}(0)=O(\varepsilon^2).
\label{eq:Energy_u1}
\end{split}
\end{equation}
The NKGE with $O(\varepsilon^p)$ nonlinearity and $O(1)$ initial data is equivalent to it with $O(1)$ nonlinearity and $O(\varepsilon)$ initial data. For simplicity, we only present the numerical method and its error estimates for the NKGE with weak nonlinearity in the following sections. Extensions of the numerical method and corresponding error estimates for the NKGE with small initial data are straightforward.

For the NKGE \eqref{eq:WNE} with $\varepsilon = 1$, there are extensive analytical and numerical results in the literature. Along the analytical front, the existence of global classical solutions and almost periodic solutions have been investigated, we refer to \cite{BP,BV,VW} and references therein. In the numerical aspect, different numerical schemes have been proposed and analyzed, including the finite difference time domain (FDTD) methods \cite{BD,CWG}, exponential wave integrator Fourier pseudospectral (EWI-FP) method \cite{BD,BDZ,BY}, fourth-order compact method \cite{DMA}, asymptotic-preserving (AP) schemes \cite{CCLM}, multiscale time integrator Fourier pseudospectral (MTI-FP) method \cite{BCZ,BZ}, etc. For comparisons of different numerical methods, we refer the readers to \cite{BZ2,JV}. However, for the NKGE \eqref{eq:WNE} with $0< \varepsilon \ll 1$, the analysis and numerical computation of the long time dynamics are mathematically rather complicated. The existence of the solution of the Cauchy problem for the NKGE with weak nonlinearity/small initial data as well as the properties of the solutions have been studied in different dimensions \cite{CHL,CHL1,KT,K}. Recently, more attentions have been devoted to analyzing the life-span of the solutions to the NKGE \eqref{eq:SIE}. The analytical results indicate that the life-span of a smooth solution to the NKGE \eqref{eq:SIE}  is at least up to the time at $O(\varepsilon^{-p})$, see, e.g., \cite{D2,DS,FZ} and references therein.

The classical error estimates are normally valid up to the time at $O(1)$. Since the life-span of the solution to the NKGE \eqref{eq:WNE} is up to the time at $O(\varepsilon^{-p})$, it is natural to establish error bounds of the numerical method for the NKGE \eqref{eq:WNE} up to the time at $O(\varepsilon^{-p})$ instead of $O(1)$. In our recent work \cite{BFY}, four explicit/semi-implicit/implicit conservative/nonconservative FDTD methods for the NKGE \eqref{eq:WNE} with a cubic nonlinearity, i.e. $p = 2$, have been proposed and anlyzed. The error estimates of the FDTD methods indicate that, in order to obtain `correct' numerical solution of the NKGE \eqref{eq:WNE} with a cubic nonlinearity, the $\varepsilon-$scalability is: $h=O(\varepsilon)$ and $\tau = O(\varepsilon)$. The fourth-order compact finite difference (4cFD) method could obtain higher order approximations with less grid points \cite{L,LZK,MAD2}, which is able to improve the resolution capacity especially for $0 < \varepsilon \ll 1$. The aim of this paper is to extend the cubic case to a general pure power case and establish rigorous error bounds of the 4cFD methods for the NKGE \eqref{eq:WNE} in the long time regime. In our error estimates, we pay particular attention to how the error bounds depend explicitly on the mesh size $h$ and time step $\tau$ as well as the small parameter $\varepsilon \in (0 ,1]$. Based on our rigorous error estimates, in order to get `correct' numerical approximations of the NKGE \eqref{eq:WNE} up to time at $O(\varepsilon^{-p})$, the $\varepsilon-$scalability of the 4cFD methods should be taken as:
\begin{equation}
h = O(\varepsilon^{p/4}) \quad \mbox{and} \quad \tau = O(\varepsilon^{p/2}), \quad p \in \mathbb{N}^+, \quad 0 < \varepsilon \leq 1,
\end{equation}
which performs better than the classical FDTD methods.

By a rescaling of time, i.e., $t \to t/\varepsilon^p$, the NKGE \eqref{eq:WNE} can be reformulated as an oscillatory NKGE whose solution propagates waves with wavelength at $O(1)$ and $O(\varepsilon^p)$ in space and time, respectively. The 4cFD methods to the NKGE \eqref{eq:WNE} and the error estimates can be extended straightforwardly to the oscillatory NKGE.

The rest of the paper is organized as follows. In Section 2, the fourth-order compact finite difference discretizations are presented for the NKGE \eqref{eq:WNE} and the properties of the stability, energy conservation and solvability are analyzed. In Section 3, we establish rigorous error bounds of the 4cFD methods for the NKGE \eqref{eq:WNE} up to the time at $O(\varepsilon^{-p})$. Numerical results are reported in Section 4 to confirm our error bounds. In Section 5, we extend the 4cFD methods and their error bounds to an oscillatory NKGE. Finally, some conclusions are drawn in Section 6. Throughout this paper, we adopt the notation $A \lesssim B$ to represent that there exists a generic constant $C > 0$, which is independent of the mesh size $h$ and time step $\tau$ as well as $\varepsilon$ such that $|A| \leq C B$.

\section{4cFD methods and their analysis}
In this section, we adapt the 4cFD methods to discretize the NKGE \eqref{eq:WNE} and analyze their stability, energy conservation and solvability. For simplicity of notations, we shall only present the numerical methods and their analysis for the NKGE (\ref{eq:WNE}) in one dimension (1D). Generalizations to higher dimensions are straightforward and results remain valid with minor modifications. In 1D, we consider the following NKGE
\begin{equation}
\begin{split}
&\partial_{tt} u(x, t) - \partial_{xx} u (x, t)+  u(x, t) + \varepsilon^p u^{p+1}(x, t) = 0,\quad x \in \Omega = (a, b),\quad t > 0, \\
&u(x, 0) =\phi(x), \quad \partial_t u(x, 0) =\gamma(x) , \quad x \in \overline{\Omega} = [a, b],
\label{eq:21}
\end{split}
\end{equation}
with periodic boundary conditions.

\subsection{4cFD methods}
Choose the time step $\tau := \Delta t > 0$ and mesh size $h : = \Delta x > 0$, and denote  $M=(b - a)/h $ being a positive integer, the grid points and time steps as
\begin{equation}
x_j := a + j h, \quad j = 0, 1, \ldots, M; \quad t_n := n\tau, \quad n = 0, 1, 2, \ldots.
\end{equation}
Define the index sets $\mathcal{T}_M = \{j | j = 0, 1,\cdots, M-1\}$, $\mathcal{T}^0_M = \{j | j = 0,1,2,\cdots, M\}$.
Denote $X_M = \{u = (u_0,u_1, \ldots, u_M)^T | u_0 = u_M\} \in \mathbb{R}^{M+1}$ and we always use $u_{-1} = u_{M-1}$ and $u_{M+1} = u_{1}$ if they are involved.
The standard discrete $l^2$, semi-$H^1$ and $l^{\infty}$ norms and inner product in $X_M$ are defined as
\begin{equation*}
\begin{split}
& \|u\|^2_{l^2} = h \sum_{j \in \mathcal{T}_M} |u_j|^2, \quad \|\delta^{+}_x u\|^2_{l^2} = h \sum_{j \in \mathcal{T}_M} |\delta^{+}_x u_j|^2,\\
& \|u\|_{l^{\infty}} = \max_{j \in \mathcal{T}_M} |u_j|, \quad  (u, v) = h \sum_{j \in \mathcal{T}_M} u_j v_j,
\end{split}
\end{equation*}
with $\delta^{+}_x u\in X_M$ defined as $\delta^{+}_x u_j=(u_{j+1}-u_j)/h$ for $j \in \mathcal{T}_M$.

Let $u^n_j$ be the numerical approximation of $u(x_j, t_n)$ for $j \in \mathcal{T}^0_M$, $n \geq 0$ and denote the numerical solution at time $t = t_n$ as $u^n=(u_0^n, u_1^n, \ldots, u_M^n)^T\in X_M$. We introduce the finite difference operators as
\begin{equation*}
\delta^{+}_t u^n_j = \frac{u^{n+1}_j-u^n_j}{\tau},\quad \delta^{-}_t u^n_j = \frac{u^{n}_j-u^{n-1}_j}{\tau},\quad \delta^{2}_t u^n_j = \frac{u^{n+1}_j-2u^{n}_j+u^{n-1}_j}{{\tau}^2},
\end{equation*}
\begin{equation*}
\delta^{+}_x u^n_j = \frac{u^{n}_{j+1}-u^n_j}{h},\quad \delta^{-}_x u^n_j = \frac{u^{n}_j-u^{n}_{j-1}}{h},\quad \delta^{2}_x u^n_j = \frac{u^{n}_{j+1}-2u^{n}_j+u^{n}_{j-1}}{h^2},
\end{equation*}
\begin{equation*}
\mathcal{A} u^n_j = \left(1+\frac{h^2}{12}\delta^2_x\right) u^n_j = \frac{1}{12}(u^n_{j-1}+10u^n_j+u^n_{j+1}) .
\end{equation*}
To simplify notations, for a function $u(x,t)$ and a grid function $u^n \in X_M (n \geq 0)$, we denote for $n \geq 1$,
\begin{equation*}
u(x,t_{[n]}) = \frac{u(x, t_{n+1})+ u(x, t_{n-1})}{2}, \ x \in \bar{\Omega};\quad u^{[n]}_j = \frac{u^{n+1}_j + u^{n-1}_j}{2}, \ j \in \mathcal{T}^0_M.
\end{equation*}

A fourth-order approximation is implemented by replacing the central difference operator $\delta^2_x$ with $(1-\frac{h^2}{12}\delta^2_x)\delta^2_x$, which requires a five-point stencil. In order to obtain a compact three-point stencil, $(1-\frac{h^2}{12}\delta^2_x)\delta^2_x$ is approximated by $(1+\frac{h^2}{12}\delta^2_x)^{-1}\delta^2_x$ \cite{L,LZK,LLG,MAD,WGX}.

In this paper, we consider the following fourth-order compact finite difference (4cFD) methods:

I. The implicit 4cFD method
\begin{equation}
\delta^2_t u^n_j - \mathcal{A}^{-1}\delta^2_x u^{[n]}_j + u^{[n]}_j + \varepsilon^p G\left(u^{n+1}_j, u^{n-1}_j \right)=0,\quad j \in \mathcal{T}_M,\quad n \geq 1.
\label{eq:I4th}
\end{equation}

II. The semi-implict 4cFD method
\begin{equation}
\delta^2_t u^n_j - \mathcal{A}^{-1}\delta^2_x u^{[n]}_j + u^{[n]}_j + \varepsilon^p(u^n_j)^{p+1}=0,\quad j \in \mathcal{T}_M,\quad n \geq 1.
\label{eq:SI4th}
\end{equation}
Here,
\begin{equation}
G(v, w) = \frac{F(v)-F(w)}{v-w}, \ v, w \in \mathbb{R}, \ F(v)=\int^v_0 s^{p+1} ds = \frac{v^{p+2}}{p+2},\ v \in \mathbb{R}.
\end{equation}
The initial and boundary conditions in \eqref{eq:21} are discretized as 
\begin{equation}
u^{n+1}_0 = u^{n+1}_M,\quad u^{n+1}_{-1} = u^{n+1}_{M-1}, \quad n \geq 0;\quad u^0_j = \phi(x_j),\quad j \in \mathcal{T}^0_M,
\label{eq:wib}
\end{equation}
where the first step $u^1$ is updated by the initial data and Taylor expansion as
\begin{equation}
u^1_j = \phi(x_j) + \tau \gamma(x_j) + \frac{{\tau}^2}{2} \left[ \phi''(x_j)- \phi(x_j)-\varepsilon^p \left(\phi(x_j)\right)^{p+1}\right], \ j \in \mathcal{T}^0_M.
\label{eq:w1}
\end{equation}
\begin{remark}
For the first step $u^1$, we can also replace $\phi''(x_j)$ by $\mathcal{A}^{-1}\delta^2_x\phi(x_j)$ when it is not easy to calculate $\phi''(x)$.
\end{remark}

Clearly, the above 4cFD methods are time symmetric or time reversible, i.e., they are unchanged if interchanging $n+1 \leftrightarrow n-1$ and $\tau \leftrightarrow -\tau$. The implicit 4cFD \eqref{eq:I4th} can be solved via either a direct solver or an iterative method. The semi-implicit 4cFD \eqref{eq:SI4th} can be explicitly updated in the Fourier space with $O(MlnM)$ computational cost per time step \cite{BC2,BFY}. For other semi-implicit and explicit schemes, we leave them to readers who are interested in.

\subsection{Some useful lemmas}

The operator $\mathcal{A}$ can be written as a matrix 
\begin{equation}
A
=\frac{1}{12}\begin{pmatrix}
10 &  1  & & 1\\
1  &  10  &1 & \\
    & \ddots & \ddots  & \\
1 &  & 1 & 10\\
\end{pmatrix}.
\end{equation}
It is easy to check that $A$ is an $(M+1) \times (M+1)$ positive definite matrix, then we can introduce a new discrete norm $\|u\|_{\ast} = \sqrt{(A^{-1}u, u)}$ for $u \in X_M$. The following lemmas will be used in our error estimates. The proof proceeds in the analogous lines as in \cite{LLG,WGX} and we just show the proof of Lemma \ref{lemma:norm_eq} in detail here for brevity.

\begin{lemma}
For any two grid functions $u, v \in X_M$, it holds
\begin{equation}
(\delta^+_x\delta^-_x u, v) = -(\delta^+_x u, \delta^+_x v).
\end{equation}
\label{lemma:delta}
\end{lemma}

\begin{lemma}
The operators $\mathcal{A}$ and $\mathcal{A}^{-1}$ are commutative with $\delta^+_x$ and $\delta^-_x$, i.e. for any grid function $u \in X_M$,
\begin{equation}
\begin{split}
&\delta^+_x \mathcal{A} u =  \mathcal{A} \delta^+_x u, \quad \delta^-_x \mathcal{A} u =  \mathcal{A} \delta^-_x u,\\
&\delta^+_x \mathcal{A}^{-1} u =  \mathcal{A}^{-1} \delta^+_x u, \quad \delta^-_x \mathcal{A}^{-1} u =  \mathcal{A}^{-1} \delta^-_x u.\\
\end{split}
\end{equation}
\label{lemma:com}
\end{lemma}

\begin{lemma}
The discrete norms $\|\cdot\|_{\ast}$ and $\|\cdot\|_{l^2}$ are equivalent. In fact, for any grid function $u \in X_M$, it holds
\begin{equation}
 \|u\|_{l^2} \leq \|u\|_{\ast} \leq \frac{\sqrt{6}}{2}\|u\|_{l^2}.
\end{equation}
\label{lemma:norm_eq}
\end{lemma}
\emph{Proof.}
For $\forall x \in \mathbb{R}^{M+1}$, $x = (x_0, x_2, \cdots, x_{M})^T $, we have
\begin{equation}
\begin{split}	
12x^T A x & = 10 \sum_{j=0}^M x_j^2 + 2 \sum^{M}_{j=1}x_{j-1}x_j +2 x_{0}x_{M}\\	
& = 10 \sum_{j=0}^M x_j^2 + \sum^{M}_{j=1} (x_{j-1} + x_j)^2 - \sum^{M}_{j=1} (x_{j-1}^2 + x_j^2) + (x_0+x_{M})^2 - (x_0^2+x_{M}^2)\\
& = 8 \sum^{M}_{j=0} x_j^2  + \sum^{M}_{j=1} (x_{j-1} + x_j)^2 + (x_0+x_{M})^2\\
& \geq 8x^Tx 
\end{split}
\label{eq:A1}
\end{equation}
and 
\begin{equation}
\begin{split}
12x^T A x & = 10 \sum^{M}_{j=0} x_j^2 + 2 \sum^{M}_{j=1}x_{j-1}x_j +2 x_{0}x_{M}\\	
& \leq 10 \sum^{M}_{j=0} x_j^2  + \sum^{M}_{j=1} (x_{j-1}^2+x_{j}^2)+ (x_0^2+x_{M}^2)\\
& = 12\sum^{M}_{j=0} x_j^2 \\
& = 12x^T x.
\end{split}
\label{eq:A2}
\end{equation}
As we know,
\begin{equation}
\lambda_{\min} x^Tx \leq x^T A x \leq	\lambda_{\max} x^T x,
\label{eq:egi_A}
\end{equation}
where $\lambda_{\min}$ and $\lambda_{\max}$ are the minimal and maximal eigenvalues of the matrix $A$, respectively. We take the left (right) equal if and only if $x$ is the eigenvector of $\lambda_{\min} (\lambda_{\max})$. From \eqref{eq:A1} and \eqref{eq:A2}, we could obtain 
 \begin{equation*}
 \frac{2}{3}\leq \lambda_{\min}\leq \lambda_{\max} \leq 1.		
 \end{equation*}
Applying \eqref{eq:egi_A} to the matrix $A^{-1}$, we have 
\begin{equation*}
\|u\|^2_{l^2} \leq \min_{\lambda_j \in \sigma(A^{-1})} \lambda_j \|u\|^2_{l^2} \leq \|u\|^2_{\ast} \leq 	\max_{\lambda_j \in \sigma(A^{-1})} \lambda_j \|u\|^2_{l^2} \leq \frac{3}{2} \|u\|^2_{l^2}
\end{equation*}
which replies 
\begin{equation*}
 \|u\|_{l^2} \leq \|u\|_{\ast} \leq \frac{\sqrt{6}}{2}\|u\|_{l^2}.
\end{equation*}
\hfill $\square$ \bigskip

\subsection{Stability, energy conservation and solvability}
Let $T_0 > 0$ be a fixed constant, and denote 
\begin{equation}
\sigma_{\max} := \max_{0 \leq n \leq T_0 \varepsilon^{-p}/\tau} \|u^n\|^p_{l^{\infty}}.
\end{equation}
By using the standard von Neumann analysis \cite{BD,BFY}, we have the following lemma for the stability of the 4cFD methods for the NKGE \eqref{eq:21}.

\begin{lemma}
(stability) For the above 4cFD methods applied to the NKGE \eqref{eq:21} up to the time $t=T_0/\varepsilon^{p}$, we have:

(i)The implicit 4cFD \eqref{eq:I4th} is unconditionally stable for any $h > 0$, $\tau > 0$ and $0 < \varepsilon \leq 1$.

(ii) When $\sigma_{\max} \leq \varepsilon^{-p}$, the semi-implicit 4cFD \eqref{eq:SI4th} is unconditionally stable for any $h > 0$ and $\tau >0$; and when $\sigma_{\max} > \varepsilon^{-p}$, this scheme is conditionally stable under the stability condition 
\begin{equation}
0 < \tau < \frac{2}{\sqrt{\varepsilon^p\sigma_{\max}-1}},\quad h > 0, \quad 0 < \varepsilon \leq 1. 
\label{eq:stability_semi}
\end{equation}
\label{lemma:stability}
\end{lemma}
\emph{Proof.} Replacing the nonlinear term by $f(u) = \varepsilon^p\sigma_{\max} u$, plugging
\begin{equation*}
u^{n-1}_j = \sum_{l} \hat{U}_l e^{2ijl\pi/M}, \quad u^{n}_j = \sum_{l} \xi_l \hat{U}_l e^{2ijl\pi/M},\quad  u^{n+1}_j = \sum_{l} \xi^2_l \hat{U}_l e^{2ijl\pi/M},
\end{equation*}
into \eqref{eq:I4th} and  \eqref{eq:SI4th}, with $\xi_l$ the amplification factor of the $l$th mode in phase space, we have the characteristic equation with the following structure
\begin{equation}
\xi^2_l - 2\theta_l \xi_l +1  = 0, \quad l = -\frac{M}{2}, \cdots, \frac{M}{2}-1,	
\label{eq:cha}
\end{equation}
where $\theta_l \in \mathbb{R}$ is determined by the 4cFD methods \eqref{eq:I4th} and  \eqref{eq:SI4th}. Solving the characteristic equation \eqref{eq:cha}, we have $\xi_l = \theta_l \pm \sqrt{\theta^2_l-1}$. The stability of numerical schemes amounts to 
\begin{equation}
|\xi_l| \leq 1 \ \Longleftrightarrow	\ |\theta_l|\leq 1,\quad  l = -\frac{M}{2}, \cdots, \frac{M}{2}-1.
\end{equation} 

(i) For the implicit 4cFD \eqref{eq:I4th}, we have 
\begin{equation}
0 \leq \theta_l = \frac{2}{2+\tau^2(1+\varepsilon^p\sigma_{\max}+c\lambda_l^2)} \leq 1, \quad l = -\frac{M}{2}, \cdots, \frac{M}{2}-1,	
\end{equation}
with
\begin{equation}
 c = \frac{3}{3-\sin^2\left(\frac{\pi l}{M}\right)} ,\quad \lambda_l = \frac{2}{h} \sin\left(\frac{\pi l}{M}\right),\quad \mu_l = \frac{2\pi l}{b-a}, \quad  l = -\frac{M}{2}, \cdots, \frac{M}{2}-1.	
\end{equation}
This implies that the implicit 4cFD \eqref{eq:I4th} is unconditionally stable for any $h > 0$, $\tau > 0$ and $0 < \varepsilon \leq 1$.

(ii) For the semi-implicit 4cFD \eqref{eq:SI4th}, we have
\begin{equation}
\theta_l = \frac{2-\tau^2\varepsilon^p \sigma_{\max}}{2+\tau^2(1+ c\lambda_l^2)},\quad l = -\frac{M}{2}, \cdots, \frac{M}{2}-1.
\end{equation}
Noticing $c \geq 1$ and $0 \leq \lambda_l^2 \leq \frac{4}{h^2}$, when $\sigma_{\max} \leq \varepsilon^{-p}$, or $\sigma_{\max} > \varepsilon^{-p}$ with the condition \eqref{eq:stability_semi}, we can get
\begin{equation*}
	\tau^2(\varepsilon^p\sigma_{\max} - 1 - c \lambda_l^2) \leq \tau^2(\varepsilon^p\sigma_{\max} - 1) \leq 4 \implies |\theta_l| \leq 1,\quad l = -\frac{M}{2}, \cdots, \frac{M}{2}-1.
\end{equation*}
The proof is completed.
\hfill $\square$ \bigskip

\begin{remark}
The stability of the semi-implicit 4cFD \eqref{eq:SI4th} is related to $\sigma_{\max}$, dependent on the boundedness of the $l^{\infty}$ norm of the numerical solution $u^n$. The error estimates up to the previous time step could ensure the boundedness, by the inverse inequality, and such an error estimate could be recovered at the next time step, as given by the Theorem presented in Section 3.
\end{remark}

The implicit 4cFD \eqref{eq:I4th} conserves the energy in the discrete level, while the semi-implicit 4cFD \eqref{eq:SI4th} does not. We have the following lemma for the energy conservation.
\begin{lemma}
(energy conservation) For $n \geq 0$, the implicit 4cFD \eqref{eq:I4th} conserves the discrete energy as 
\begin{equation}
\begin{split}
E^n := &\|\delta^+_t u^n\|^2_{l^2} + \frac{1}{2} \left( \|\delta^+_x u^n\|^2_{\ast} + \|\delta^+_x u^{n+1}\|^2_{\ast}\right) + \frac{1}{2} \left( \|u^n\|^2_{l^2} + \|u^{n+1}\|^2_{l^2}\right) \\
& + \frac{\varepsilon^p h}{p+2} \sum_{j \in \mathcal{T_M}} \left[|u^n_j|^{p+2}+|u^{n+1}_j|^{p+2}\right] \equiv E^0.
\end{split}
\end{equation}
\label{lemma:conservation}
\end{lemma}

\begin{lemma}
(solvability of the 4cFD methods) For any given $u^n$, $u^{n-1} (n \geq 1)$, there exists a unique solution $u^{n+1}$ of the 4cFD methods \eqref{eq:I4th} and \eqref{eq:SI4th} with \eqref{eq:wib} and \eqref{eq:w1}.
\label{lemma:solvability}
\end{lemma}

The proofs of Lemmas \ref{lemma:conservation} and \ref{lemma:solvability} proceed in the analogous lines as in \cite{BD,BFY,BS1} and we omit the details here for brevity.

\section{Error estimates for the 4cFD methods}
In this section, we will rigorously establish the error bounds of the 4cFD methods for the NKGE \eqref{eq:21}.

\subsection{Main results}
According to the known results in \cite{D2,DS,K} and references therein, we can make the assumptions on the exact solution $u$ of the NKGE \eqref{eq:21} up to the time $t = T_0/\varepsilon^p$:
\begin{equation*}
(A)
\begin{split}
u \in \ &  C([0, T_0/\varepsilon^p]; W_p^{6, \infty})  \cap C^2([0, T_0/\varepsilon^p]; W^{4, \infty}) \cap C^4([0, T_0/\varepsilon^p]; W^{2, \infty}), \\ &\quad\left\|\frac{\partial^{r+q}}{\partial t^r \partial x^q} u(x, t)\right\|_{L^{\infty}}  \lesssim 1,\quad 0 \leq r \leq 4,\quad 0 \leq r+q \leq 6,
\end{split}
\end{equation*}
here $W^{m, \infty}_p = \{ u \in W^{m, \infty}| \frac{\partial^l}{\partial x^l} u(a) = \frac{\partial^l}{\partial x^l} u(b),\quad 0\leq l< m \}$ for $m \geq 1$.

Denote  $M_0 = \sup_{\varepsilon\in(0,1]} \| u(x, t)\|_{L^{\infty}}$ and the error function ${e}^n \in X_M (n \geq 0)$ as
\begin{equation}
e^n_j = u(x_j, t_n) - u^n_j, \quad j \in \mathcal{T}^0_M, \quad n \geq 0, 
\label{eq:verrdef}
\end{equation}
where $u^n\in X_M$ is the numerical approximation of the NKGE (\ref{eq:21}), then we have the following error estimates for the implicit 4cFD \eqref{eq:I4th} with \eqref{eq:wib} and \eqref{eq:w1}:

\begin{theorem}
Under the assumption ($A$), there exist constants $h_0 > 0$ and $\tau_0 > 0$ sufficiently small and independent of $\varepsilon$, such that for any $0 < \varepsilon \leq 1$, when  $0 < h \leq h_0\varepsilon^{p/4}$, $0 < \tau \leq \tau_0 \varepsilon^{p/2}$, the following two error estimates of the scheme (\ref{eq:I4th}) with (\ref{eq:wib}) and (\ref{eq:w1}) hold 
\begin{equation}
\|e^n\|_{l^2} +\|\delta^{+}_x e^n\|_{l^2} \lesssim \frac{h^4}{\varepsilon^p}+ \frac{\tau^2}{\varepsilon^p},\quad \|u^n\|_{l^{\infty}} \leq 1 + M_0, \quad 0 \leq n \leq \frac{T_0/\varepsilon^p}{\tau}.
\label{eq:eb}
\end{equation}
\label{thm:I4th}
\end{theorem}

For the semi-implict 4cFD \eqref{eq:SI4th} with \eqref{eq:wib} and \eqref{eq:w1}, the error bounds can be established as follows:
\begin{theorem}
Under the assumption ($A$), there exist constants $h_0 > 0$ and $\tau_0 > 0$ sufficiently small and independent of $\varepsilon$, such that for any $0 < \varepsilon \leq 1$, when  $0 < h \leq h_0\varepsilon^{p/4}$, $0 < \tau \leq \tau_0 \varepsilon^{p/2}$ and under the stability condition \eqref{eq:stability_semi}, the following two error estimates of the scheme (\ref{eq:SI4th}) with (\ref{eq:wib}) and (\ref{eq:w1}) hold 
\begin{equation}
\|e^n\|_{l^2} +\|\delta^{+}_x e^n\|_{l^2} \lesssim \frac{h^4}{\varepsilon^p}+ \frac{\tau^2}{\varepsilon^p},\quad \|u^n\|_{l^{\infty}} \leq 1 + M_0, \quad 0 \leq n \leq \frac{T_0/\varepsilon^p}{\tau}.
\label{eq:eb}
\end{equation}
\label{thm:SI4th}
\end{theorem}

\begin{remark}
The above error bounds in Theorem \ref{thm:I4th} and Theorem \ref{thm:SI4th} are still valid in higher dimensions, e.g., $d = 2, 3$, provided that the technical conditions $0 < h \lesssim \varepsilon^{p/4}\sqrt{C_d(h)}$ and $0 < \tau \lesssim \varepsilon^{p/2} \sqrt{C_d(h)}$ with
\begin{equation}
C_d(h)=
\begin{cases}
1/|\ln h|, &d=2,\\
h^{1/2}, & d=3.
\end{cases}
\end{equation}
The reason is due to the discrete Sobolev inequality \cite{BC2,BD,BS,T}:
\begin{equation}
 \|u^n\|_{l^{\infty}} \lesssim \frac{1}{C_d(h)} \left( \|u^n\|_{l^2} + \|\delta^+_x u^n\|_{l^2}\right).
 \end{equation}
\end{remark}

Based on the above theorems, the 4cFD methods have the following spatial/temporal resolution capacity for the NKGE \eqref{eq:21} in the long time regime. In fact,  given an accuracy bound $\delta_0 > 0$, the $\varepsilon$-scalability of the 4cFD methods is:
\begin{align}
h = O(\varepsilon^{p/4}\sqrt{\delta_0}) =O(\varepsilon^{p/4}), \quad \tau = O(\varepsilon^{p/2} \sqrt{\delta_0}) =O(\varepsilon^{p/2}), \quad  0 < \varepsilon \leq 1.
\label{eq:WNE_scalability}
\end{align}

Compared with the commonly used standard FDTD methods \cite{BFY}, the results can attain higher order accuracy in space for a given mesh size or improve the spatial resolution capacity, i.e., it needs less grid points while maintaining the same accuracy.

\subsection{Proof of Theorem \ref{thm:SI4th}}
For the semi-implicit 4cFD method, we establish the error bounds in Theorem \ref{thm:SI4th} by the method of mathematical induction \cite{BC2,BCJY}. Throughout this section, the stability condition \eqref{eq:stability_semi} is assumed. 

Denote the local truncation error as $\xi^n \in X_M$ for $0 \leq n \leq \frac{T_0/\varepsilon^p}{\tau} - 1$
\begin{equation}
\begin{split}
{\xi}^0_j :=& \ \delta^{+}_t u(x_j, 0) - \gamma(x_j) - \frac{\tau}{2}\left[\phi''(x_j) -   \phi(x_j) - \varepsilon^p (\phi(x_j))^{p+1}\right],\quad j \in \mathcal{T}_M,\\
{\xi}^n_j := &\  \delta^{2}_t u(x_j, t_n) - \mathcal{A}^{-1} \delta^2_x u(x_j, t_{[n]})+  u(x_j, t_{[n]}) + \varepsilon^pu^{p+1}(x_j, t_n),\quad n \geq 1,
\end{split}
\label{eq:WNE_lt}
\end{equation}
and the error of the nonlinear term as $\eta^n \in X_M$
\begin{equation}
\eta^n_j := \varepsilon^p\left(u^{p+1}(x_j, t_n) - (u^n_j)^{p+1}\right),\quad j \in \mathcal{T}_M, \quad 1 \leq n \leq \frac{T_0/\varepsilon^p}{\tau} - 1.
\label{eq:nonlinear_def}
\end{equation}

We begin with the error estimates of local truncation error ${\xi}^n \in X_M$.
\begin{lemma}
Under the assumption (A), we have
\begin{equation}
\|\xi^0\|_{l^2} + \|\delta^+_x \xi^0\|_{l^2} \lesssim  \tau^2,\quad \|\xi^n\|_{l^2} \lesssim h^4+\tau^2,\quad 1 \leq n \leq \frac{T_0/\varepsilon^p}{\tau}-1.
\label{eq:lterr}
\end{equation}
\label{thm:Iterr}
\end{lemma}
\emph{Proof.}
Under the assumption (A), by applying the Taylor expansion, Young's inequality and Lemma \ref{lemma:com}, we have 
\begin{equation*}
\begin{split}
|{\xi}^0_j| \lesssim & \ \tau^2 \|\partial_{ttt} u \|_{L^{\infty}} \lesssim \tau^2, \quad j \in \mathcal{T}_M, \\
|\mathcal{A}{\xi}^n_j| \lesssim & \ \tau^2 \left[ \| \partial_{tt} u \|_{L^{\infty}} + \| \partial_{tttt} u \|_{L^{\infty}}  + \| \partial_{ttxx} u \|_{L^{\infty}} + \| \partial_{ttttxx} u \|_{L^{\infty}} \right] \\
&+ h^4 \left[\| \partial_{xxxx} u \|_{L^{\infty}} + \| \partial_{ttxxxx} u \|_{L^{\infty}} + \| \partial_{xxxxxx} u \|_{L^{\infty}} \right]\lesssim \ h^4+\tau^2, \ n \geq 1,
\end{split}
\end{equation*}
which leads to $|{\xi}^n_j| \lesssim \ h^4+\tau^2$ for $j \in \mathcal{T_M}, \ n \geq 1$. Similarly, we have $|\delta^{+}_x {\xi}^0_j| \lesssim \tau^2$ for $j \in \mathcal{T_M}$. These immediately imply \eqref{eq:lterr}.
\hfill $\square$ \bigskip

Next, the error equation for the semi-implicit 4cFD \eqref{eq:SI4th} is
\begin{equation}
\begin{split}
&\delta^2_t e^n_j - \mathcal{A}^{-1} \delta^2_x e^{[n]}_j + e^{[n]}_j = \xi^n_j - \eta^n_j,\quad 1 \leq n \leq \frac{T_0/\varepsilon^p}{\tau},\\
&e^0_j = 0,\quad e^1_j=\tau \xi^0_j, \quad j \in \mathcal{T}_M.
\end{split}
\label{eq:errorfun}
\end{equation}

We will improve Theorem \ref{thm:SI4th} by the method of mathematical induction. For $n = 0$, \eqref{eq:eb} is trivial. For $n = 1$, the error function \eqref{eq:errorfun} and the error estimates of the local truncation error \eqref{eq:lterr} imply 
\begin{equation*}
\|e^1\|_{l^2}=\tau\|\xi^{0}\|_{l^2} \leq C_1 \tau^3, \quad \|\delta_{x}^{+} e^1\|_{l^2}=\tau \|\delta_{x}^{+} {\xi}^{0}\|_{l^2} \leq C_1 \tau^3.
\end{equation*}
By the triangle inequality, discrete Sobolev inequality and the assumption (A), there exist a constant $\tau_1 > 0$ sufficiently small, when $0 < \tau < \tau_1$, we have 
\begin{equation*}
\|u^{1}\|_{l^{\infty}} \leq \|u\left(x, t_{1}\right)\|_{L^{\infty}}+ \|e^{1}\|_{l^{\infty}} \leq \|u\left(x, t_{1}\right)\|_{L^{\infty}}+\|e^{1}\|_{l^{2}}+\|\delta_{x}^{+} e^{1}\|_{l^2} \leq M_{0}+1,
\end{equation*}
which immediately implies \eqref{eq:eb} for $n = 1$.

Now assuming that \eqref{eq:eb} is valid for all $0 \leq n \leq m-1 \leq \frac{T_0/\varepsilon^p}{\tau}-1$, it needs to prove that it is still valid when $n = m$. Under the assumption (A), the error of the nonlinear term $\eta^n$ for $1 \leq n \leq m-1$ can be controlled as
\begin{equation}
 \|\eta^n\|_{l^2} \leq C_2 \varepsilon^p\|e^n\|_{l^2}, \quad 1 \leq n \leq m-1.
 \label{eq:eta}
\end{equation}
Define the `energy' for the error vector $e^n \in X_M (n \geq 0)$ as
\begin{equation*}
S^n = \| \delta^+_t e^n\|^2_{l^2} + \frac{1}{2} \left( \|\delta^+_x e^n\|^2_{\ast} + \|\delta^+_x e^{n+1}\|^2_{\ast}\right) + \frac{1}{2}\left(\|e^n\|^2_{l^2} + \|e^{n+1}\|^2_{l^2}\right), \quad n \geq 0.
\end{equation*}
By Lemma \ref{lemma:norm_eq}, it is easy to see that
\begin{equation*}
S^0 = \|\delta^+_t e^0\|^2_{l^2} + \frac{1}{2} \|\delta^+_x e^1\|^2_{\ast} + \frac{1}{2} \|e^1\|^2_{l^2} \lesssim  \tau^4.
\end{equation*}
Multiplying both sides of \eqref{eq:errorfun} by $h\left(e^{n+1}_j - e^{n-1}_j\right)$, summing up for $j$ and using the Young's inequality, the inequality \eqref{eq:eta}, the Lemma \ref{lemma:delta}, \ref{lemma:com} and \ref{thm:Iterr}, we derive
\begin{equation}
\begin{split}
S^n - S^{n-1} & = h \sum^{M-1}_{j=0} \left(\xi^n_j-\eta^n_j\right) \left(e^{n+1}_j - e^{n-1}_j\right) \\
& \leq \tau \varepsilon^{-p} \left(\|\xi^n\|^2_{l^2} + \|\eta^n\|^2_{l^2}\right) + \tau \varepsilon^p\left(\|\delta^+_t e^n\|^2_{l^2} + \|\delta^+_t e^{n-1}\|^2_{l^2}\right) \\
& \lesssim \tau \varepsilon^p \left(S^n+S^{n-1}\right) + \tau\varepsilon^{-p}\left(h^4 + \tau^2 \right)^2,\quad 1 \leq n \leq m-1. 
\end{split}
\end{equation}
Summing the above inequalities from 1 to $m-1$, there exists a constant $C_3 > 0$ such that
\begin{equation}
S^{m-1} \leq S^0 + C_3 \tau \varepsilon^p \sum^{m-1}_{n=0} S^n + C_3 T_0\varepsilon^{-2p} \left(h^4 + \tau^2 \right)^2. 
\end{equation}
Then the discrete Gronwall's inequality \cite{G} suggests that there exists a constant $\tau_2 > 0$ sufficiently small such that when $0 < \tau \leq \tau_2$, the following holds
\begin{equation}
S^{m-1}  \leq \left(S^0 + C_3{T_0}{\varepsilon^{-2p}}\left(h^4 + \tau^2 \right)^2 \right)e^{2C_3m\varepsilon^p\tau} \lesssim  {\varepsilon^{-2p}}\left(h^4 + \tau^2 \right)^2.
\label{eq:Sn}
\end{equation}
From the definition of $S^{m-1}$ and Lemma \ref{lemma:norm_eq}, we can obtain that $ \|e^{m}\|_{l^2}^2 +\|\delta^{+}_x e^{m}\|_{l^2}^2 \leq C_4 S^{m-1}$ when $0 < \varepsilon \leq 1$, which immediately implies
\begin{equation}
 \|e^{m}\|_{l^2} +\|\delta^{+}_x e^{m}\|_{l^2} \lesssim  \frac{h^4}{\varepsilon^p} + \frac{\tau^2}{\varepsilon^p}.
\end{equation}
The first inequality in \eqref{eq:eb} is valid for $n = m$ and it remains to estimate $\|u^m\|_{l^{\infty}}$. In fact, the discrete Sobolev inequality implies
\begin{equation}
\|e^m\|_{l^{\infty}} \lesssim \|e^m\|_{l^2} + \|\delta^+_x e^m\|_{l^2} \lesssim \frac{h^4}{\varepsilon^p} + \frac{\tau^2}{\varepsilon^p}.
\end{equation}
Thus, there exist $h_0 > 0$ and $\tau_3 > 0$ sufficiently small, when $0 < h \leq h_0 \varepsilon^{p/4} $ and  $0 < \tau \leq \tau_3 \varepsilon^{p/2} $, we obtain
\begin{equation}
\|u^m\|_{l^{\infty}} \leq \|u(x, t_m)\|_{L^{\infty}} + \|e^m\|_{l^{\infty}}  \leq  M_0+1, \quad 1 \leq m \leq T_0\varepsilon^{-p}/\tau.
\end{equation}
Under the stability condition \eqref{eq:stability_semi} and the choice of  $\tau_0 = \min\{\tau_1, \tau_2, \tau_3\}$, the estimates in \eqref{eq:eb} are valid when $n=m$. Hence, the proof of Theorem \ref{thm:SI4th} is completed by the method of mathematical induction.

\begin{remark}
For the proof of Theorem \ref{thm:I4th}, we just give the outline here and omit the details the for brevity. The key of the proof is to use the cut-off technique to deal with the nonlinearity and overcome the difficulty in uniformly bounding the numerical solution \cite{BFY,BS1,BS}. Firstly, we truncate the nonlinearity by a global Lipschitz function with compact support. Secondly, using the analogous energy method, we establish the error bounds if the exact solution is bounded and the numerical solution is close to it under some conditions. Finally, we obtain the error bounds for the implicit 4cFD method by the solvability and uniqueness of the scheme.
\end{remark}


\section{Numerical results}
In this section,  we present numerical results for the NKGE (\ref{eq:21}) up to the long time at $O(\varepsilon^{-p})$ by our proposed 4cFD methods. In the numerical experiments, we take $p=2$, $a=0$, $b=2\pi$ and choose the initial data as
\begin{align}\label{itTd1}
\phi(x) =3/(2+\cos(x)),\quad \quad\gamma(x) = \sin(x), \qquad 0\le x\le 2\pi.
\end{align}
Denote $u^n_{h,\tau}$ as the numerical solution at time $t=t_n$ obtained by the semi-implicit 4cFD method with mesh size $h$ and time step $\tau$. To quantify the numerical errors, we introduce the error function as follows:
\begin{equation}
e_{h,\tau}(t_n) = \sqrt{\|u(\cdot, t_n) - u^n_{h,\tau}\|^2_{l^2}+ \|\delta^{+}_x (u(\cdot, t_n) - u^n_{h,\tau})\|^2_{l^2}}.
\end{equation}

The `exact' solution is obtained numerically by the exponential wave integrator Fourier pseudospectral method \cite{BD,DXZ} with a very fine mesh size and a very small time step, e.g. $h_e=\pi/256$ and $\tau_e = 2\times10^{-5}$. The errors are displayed at $t = 1/\varepsilon^2$. For spatial error analysis, the time step is set as $\tau = 2 \times 10^{-5}$ such that the temporal error can be neglected; for temporal error analysis, we set the mesh size as $h = \pi/256$ such that the spatial error can be ignored.

Tables \ref{tab:WNE_h} and \ref{tab:WNE_t} present the spatial and temporal errors for different  $0 < \varepsilon \leq 1$, respectively. From Tables \ref{tab:WNE_h} and \ref{tab:WNE_t} and additional similar numerical results not shown here for brevity, we can draw the following observations:

(i) For any fixed $\varepsilon = \varepsilon_0 > 0$, the 4cFD methods are fourth-order accurate in space and second-order accurate in time (cf. the first rows in Tables  \ref{tab:WNE_h} and \ref{tab:WNE_t}). (ii) In the long time regime, the fourth order convergence in space and second order convergence in time can be observed only when $0 < h \lesssim \varepsilon^{p/4}$ and $0 < \tau \lesssim \varepsilon^{p/2}$ (cf. upper triangles above the diagonals (corresponding to $h \sim \varepsilon^{p/4}$ and $\tau \sim \varepsilon^{p/2}$, and being labelled in bold letters) in Tables  \ref{tab:WNE_h} and \ref{tab:WNE_t}), which again confirm our error estimates. In summary, our numerical results confirm our rigorous error estimates and show that they are sharp.

\begin{table}[ht!]
\renewcommand{\arraystretch}{1.5}
\caption{Spatial errors of the semi-implicit 4cFD (\ref{eq:SI4th}) for the NKGE (\ref{eq:21}) with \eqref{itTd1}}
\centering
\begin{tabular}{cccccc}
\hline
$e_{h,\tau_e}(t=1/\varepsilon^2)$ &$h_0 = \pi/8$ & $h_0/2$ &$h_0/2^2$ & $h_0/2^3$ & $h_0/2^4$ \\
\midrule[0.8pt]
$\varepsilon_0 = 1$  & 1.50E-2 & \bf{9.58E-4} & 6.02E-5 & 3.75E-6 & 2.37E-7   \\
Order & - & \bf{3.97} & 3.99 & 4.00 & 3.98  \\
$\varepsilon_0/2^2 $  & 1.02E-1 & 7.23E-3 & \bf{5.00E-4} & 3.21E-5 & 1.99E-6   \\
Order & - & 3.82 & \bf{3.85} & 3.96 & 4.01  \\
$\varepsilon _0/2^4 $  & 7.80E-1 & 6.89E-2 & 8.26E-3 & \bf{5.05E-4} & 3.17E-5   \\
Order & - & 3.50 & 3.06 & \bf{4.03} & 3.99 \\
$\varepsilon_0/2^6 $  & 5.13E-1 & 5.48E-1 & 8.43E-2 & 7.49E-3 & \bf{4.83E-4}  \\
Order & - & -0.10 & 2.70 & 3.49 & \bf{3.95}   \\
\hline
\end{tabular}
\label{tab:WNE_h}
\end{table}

\begin{table}[ht!]
\renewcommand{\arraystretch}{1.5}
\caption{Temporal errors of the semi-implicit 4cFD (\ref{eq:SI4th}) for the NKGE (\ref{eq:21}) with \eqref{itTd1}}
\centering
\begin{tabular}{ccccccc}
\hline
$e_{h_e,\tau}(t=1/\varepsilon^2)$ &$\tau_0 = 0.2 $ & $\tau_0/2 $ &$\tau_0/2^2 $ & $\tau_0/2^3 $ & $\tau_0/2^4$ & $\tau_0/2^5$ \\
\midrule[0.8pt]
$\varepsilon_0 = 1$ &\bf{6.02E-2} & 1.66E-2 & 4.32E-3 &  1.10E-3 & 2.77E-4 & 6.96E-5  \\
Order & \bf{-} & 1.86 & 1.94 & 1.97 & 1.99 & 1.99  \\
$\varepsilon_0 /2$ & 2.23E-1 & \bf{5.92E-2} & 1.50E-2 &  3.78E-3 & 9.46E-4 & 2.37E-4 \\
Order & - & \bf{1.91} & 1.98 & 1.99 & 2.00 & 2.00 \\
$\varepsilon_0 /2^2$ & 5.25E-1  & 1.46E-1 & \bf{3.86E-2} &  9.85E-3 & 2.48E-3 & 6.20E-4 \\
Order & - & 1.85 & \bf{1.92} & 1.97 & 1.99 & 2.00 \\
$\varepsilon_0 /2^3$ & 1.40E+0  & 4.82E-1 & 1.14E-1 & \bf{2.80E-2} & 7.03E-3 & 1.76E-3\\
Order & - & 1.54 & 2.08 & \bf{2.03} & 1.99 & 2.00 \\
$\varepsilon_0 /2^4$ & 3.26E+0  & 1.56E+0 & 4.83E-1 & 1.17E-1 & \bf{2.95E-2} & 7.36E-3 \\
Order & - & 1.06 & 1.69 & 2.05 &\bf{1.99} & 2.00 \\
\hline
\end{tabular}
\label{tab:WNE_t}
\end{table}

\section{Extension to an oscillatory NKGE}

\begin{figure}[ht!]
\begin{minipage}{0.5\textwidth}
\centerline{\includegraphics[width=8cm,height=6cm]{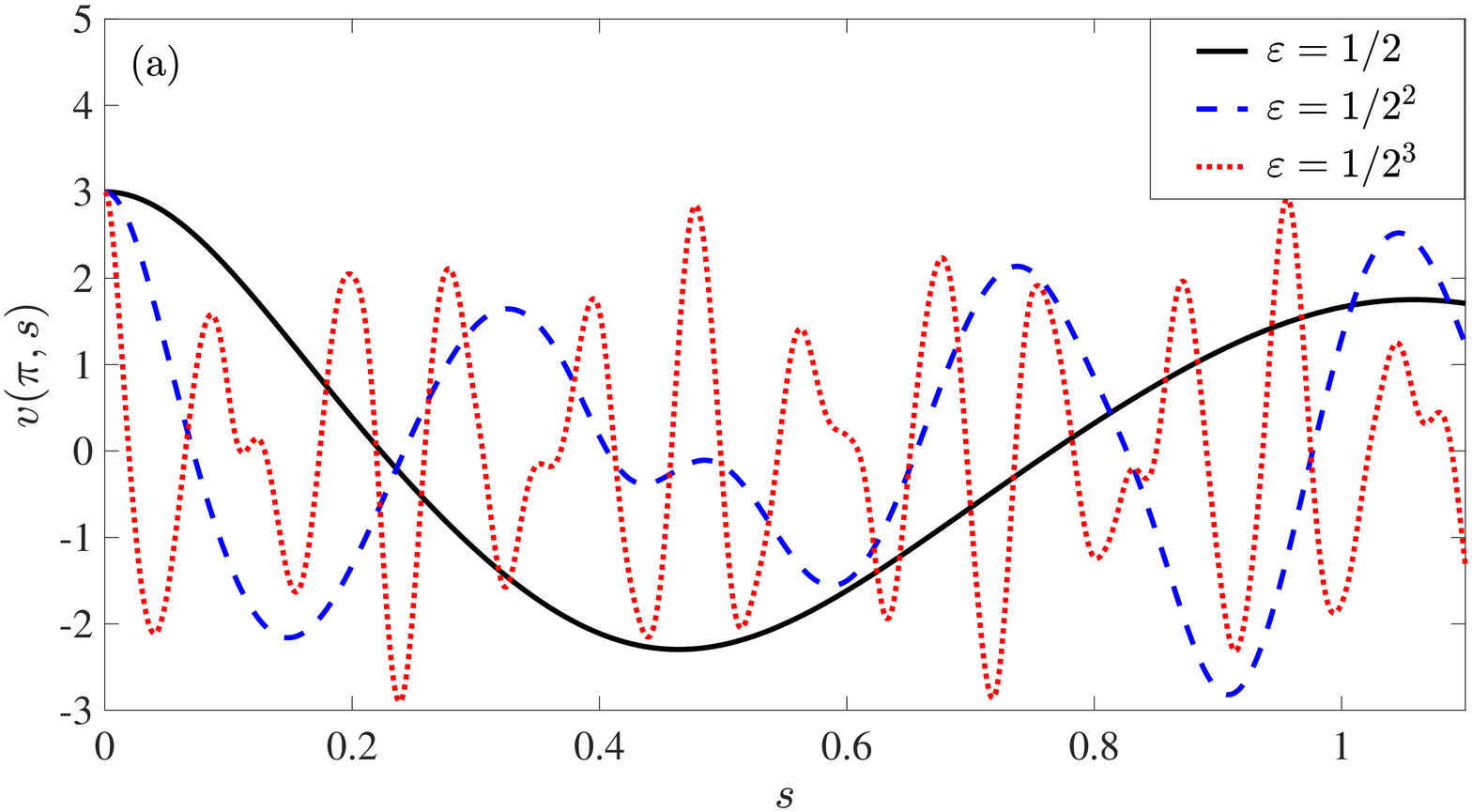}}
\end{minipage}
\begin{minipage}{0.5\textwidth}
\centerline{\includegraphics[width=8cm,height=6cm]{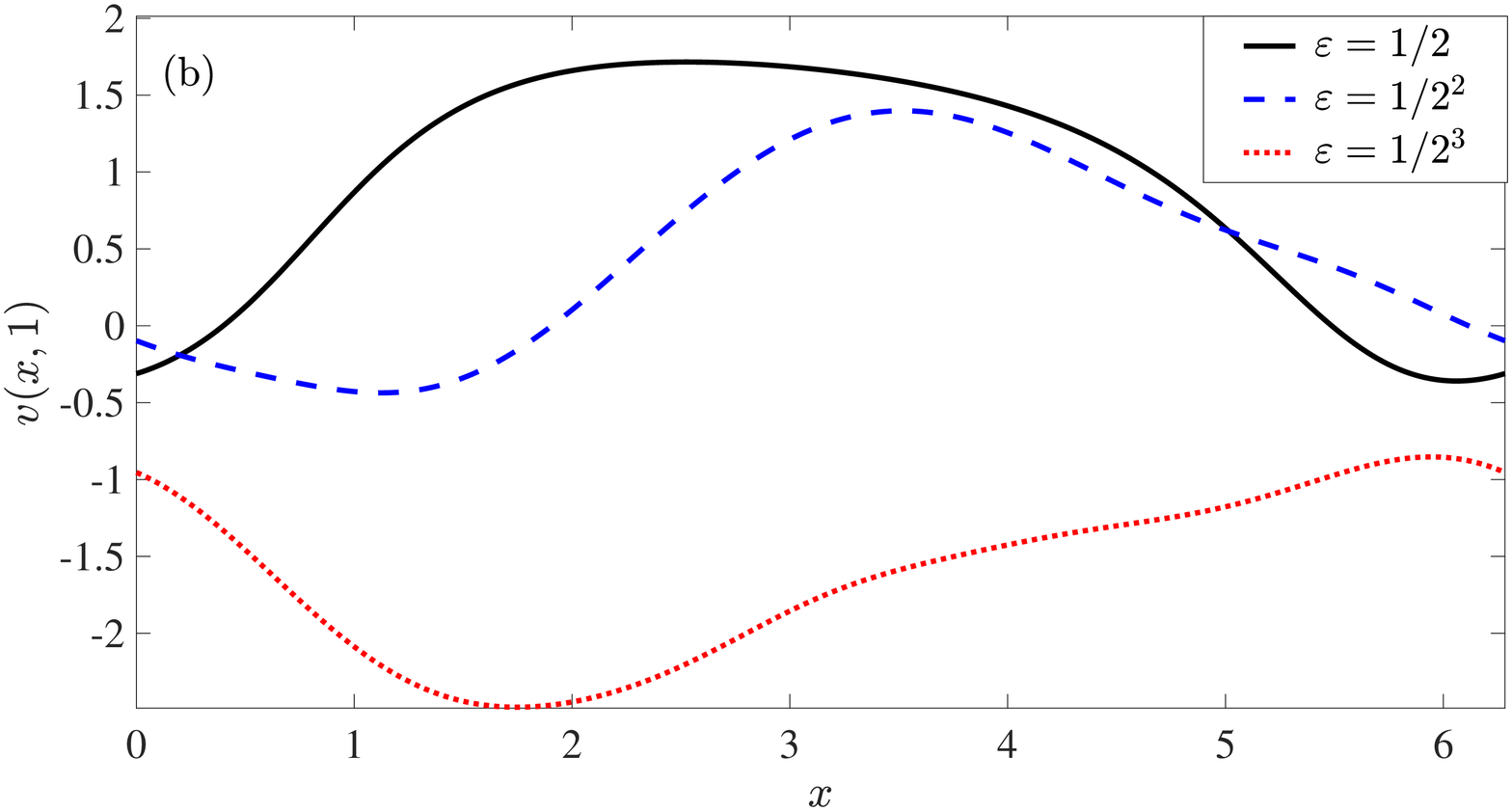}}
\end{minipage}
\vspace{-0.1in}
\caption{The solutios of the oscillatory NKGE (\ref{eq:50}) with
$d=1$, $p = 2$ and initial data \eqref{itTd1}  for different $\varepsilon$: (a) $v(\pi, s)$, (b) $v(x, 1)$.}
\label{fig:v}
\end{figure}

By a rescaling  in time $s = {\varepsilon}^p t$  and denoting $v({\bf x},s) :=u({\bf x},s/\varepsilon^p)= u({\bf x},t)$, we can reformulate the NKGE \eqref{eq:WNE}  into the following oscillatory NKGE
\begin{equation}
\begin{split}
&\varepsilon^{2p}\partial_{ss} v({\bf x,} s) - \Delta v({\bf x}, s) + v({\bf x}, s) + \varepsilon^{p} v^{p+1}({\bf x}, s) = 0,\quad {\bf{x}} \in \mathbb{T}^d,\quad s > 0,\\
&v({\bf x}, 0) = \phi({\bf x}), \quad \partial_s v({\bf x}, 0) = {\varepsilon^{-p}} \gamma({\bf x}), \quad {\bf{x}} \in \mathbb{T}^d,
\label{eq:50}
\end{split}
\end{equation}
which is also time symmetric or time reversible and conserves the {\sl energy}, i.e.,
\begin{equation}
\begin{split}
\mathcal{E}(s) := & \int_{\mathbb{T}^d} \left[ \varepsilon^{2p}|\partial_s v ({\bf{x}}, s)|^2 + |\nabla v({\bf{x}}, s)|^2 + |v({\bf{x}}, s)|^2 +\frac{2\varepsilon^p}{p+2} |v({\bf{x}}, s)|^{p+2} \right] d {\bf{x}} \\
 \equiv & \int_{\mathbb{T}^d} \left[ |\gamma({\bf{x}})|^2 + |\nabla \phi({\bf{x}})|^2 + |\phi({\bf{x}})|^2 +\frac{2\varepsilon^p}{p+2} |\phi({\bf{x}})|^{p+2}  \right] d {\bf{x}} \\
  = & E(0)=O(1), \quad s \geq 0.
\end{split}
\label{eq:Energy_v}
\end{equation}
\noindent
In fact, the long time dynamics of the NKGE \eqref{eq:WNE} up to the time at $O(\varepsilon^{-p})$ is equivalent to the dynamics of the oscillatory NKGE \eqref{eq:50}
up to the fixed time at $O(1)$. The solution of the NKGE \eqref{eq:WNE} propagates waves with wavelength at $O(1)$ in both space and time, and wave speed in space at $O(1)$ too, while the solution of the oscillatory NKGE \eqref{eq:50} propagates waves with wavelength at $O(1)$ and $O(\varepsilon^p)$ in space and time, respectively. To illustrate this, Figure \ref{fig:v} shows the solutions $v(\pi, s)$ and $v(x,1)$ of the oscillatory NKGE \eqref{eq:50} with $d=1$, $p = 2$, $\mathbb{T}=(0,2\pi)$ and initial data \eqref{itTd1} for different $0<\varepsilon\le 1$.  

In the following, we extend the 4cFD methods and their error bounds for the NKGE \eqref{eq:WNE} in Sections 2\&3 to the oscillatory NKGE (\ref{eq:50}). For simplicity of notations, the 4cFD methods and their error bounds are only presented in 1D, and the results can be generalized to higher dimensions with minor modifications. In addition, the proof of the error bounds is quite similar to that in Section 3, and the details are omitted  for brevity. We adopt similar notations as those used in Sections 2\&3 except stated otherwise. In 1D, we consider the following oscillatory NKGE
\begin{equation}
\begin{split}
&\varepsilon^{2p} \partial_{ss} v(x, s) - \partial_{xx}v(x, s) + v(x, s) + \varepsilon^{p} v^{p+1}(x, s) = 0,\ x \in \Omega = (a, b),\ s > 0,\\
&v(x, 0) = \phi(x), \quad \partial_s v(x, 0) = {\varepsilon^{-p}} \gamma(x) , \quad x \in \overline{\Omega} = [a, b],
\end{split}
\label{eq:51}
\end{equation}
with periodic boundary conditions.

\subsection{4cFD methods}
Choose the temporal step size $k := \Delta s > 0$ and denote time steps as $s_n := nk$ for $n \geq 0$. Let $v^n_j$ be the numerical approximation of $v(x_j, s_n)$ for $j \in \mathcal{T}^0_M$, $n \geq 0$ and denote the numerical solution at time $s = s_n$ as $v^n \in X_M$. Similarly, we consider the following 4cFD methods:

I. The implicit 4cFD method
\begin{equation}
\varepsilon^{2p}\delta^2_s v^n_j - \mathcal{A}^{-1} \delta^2_x v^{[n]}_j + v^{[n]}_j + \varepsilon^p G\left(v^{n+1}_j, v^{n-1}_j\right)=0,\ j \in \mathcal{T}_M, \ n \geq 1.
\label{eq:I4th_HOE}
\end{equation}

II. The semi-implicit 4cFD method
\begin{equation}
\varepsilon^{2p} \delta^2_s v^n_j - \mathcal{A}^{-1}\delta^2_x v^{[n]}_j + v^{[n]}_j + \varepsilon^p(v^n_j)^{p+1}=0,\ j \in \mathcal{T}_M, \ n \geq 1.
\label{eq:SI4th_HOE}
\end{equation}
The initial and boundary conditions in \eqref{eq:51} are discretized as 
\begin{equation}
v^{n+1}_0 = v^{n+1}_M,\quad v^{n+1}_{-1} = v^{n+1}_{M-1}, \quad n \geq 0;\quad v^0_j = \phi(x_j),\quad j \in \mathcal{T}^0_M,
\label{eq:vib}
\end{equation}
where the first step $v^1$ is updated by the Taylor expansion as
\begin{equation}
v^1_j = \phi(x_j) + \frac{k}{\varepsilon^p}\gamma(x_j) + \frac{{k}^2}{2\varepsilon^{2p}} \left[\phi''(x_j)- \phi(x_j)-\varepsilon^p \left(\phi(x_j)\right)^{p+1}\right], \ j  \in \mathcal{T}^0_M.
\label{eq:v1}
\end{equation}

We remark here that in the approximation of the first step vale, in order to uniformly bound $v^1 \in X_M$ for $\varepsilon \in (0, 1]$, $k\varepsilon^{-p}$ and $k^2\varepsilon^{-2p}$ are replaced by $\sin(k\varepsilon^{-p})$ and $k\sin(k\varepsilon^{-2p})$, respectively \cite{BFY}.

\subsection{Stability and energy conservation}
Let $T_0 > 0$ be a fixed constant, and denote 
\begin{equation}
\tilde\sigma_{\max} := \max_{0 \leq n \leq T_0/k} \|v^n\|^p_{l^{\infty}}.
\end{equation}
Similarly, we can conclude the stability of the above 4cFD methods for the oscillatory NKGE \eqref{eq:51} up to the fixed time $s=T_0$ in the following lemma.

\begin{lemma}
(stability) For the above 4cFD methods applied to the oscillatory NKGE \eqref{eq:51} up to the fixed time $s=T_0$, we have:

(i) The implicit 4cFD \eqref{eq:I4th_HOE} is unconditionally stable for any $h > 0$, $k > 0$ and $0 < \varepsilon \leq 1$.

(ii)When $\tilde\sigma_{\max} \leq \varepsilon^{-p}$, the semi-implicit 4cFD \eqref{eq:SI4th_HOE} is unconditionally stable for any $h > 0$ and $k >0$; and when $\tilde\sigma_{\max} > \varepsilon^{-p}$, this scheme is conditionally stable under the stability condition 
\begin{equation}
0 < k < \frac{{2\varepsilon^p}}{\sqrt{\varepsilon^p\tilde\sigma_{\max}-1}},\quad h > 0, \quad 0 < \varepsilon \leq 1. 
\label{eq:stability_HOE}
\end{equation}
\end{lemma}

For the implicit 4cFD \eqref{eq:I4th_HOE}, we have the following conservation property:
\begin{lemma}
The implicit 4cFD \eqref{eq:I4th_HOE} conserves the discrete energy as
\begin{equation}
\begin{split}
\mathcal{E}^n = & \varepsilon^{2p}\|\delta^+_s v^n\|^2_{l^2} + \frac{1}{2}\left(\|\delta^+_x v^n\|^2_{\ast} + \|\delta^+_x v^{n+1}\|^2_{\ast}\right) + \frac{1}{2}\left(\|v^n\|^2_{l^2} + \| v^{n+1}\|^2_{l^2}\right) \\
& + \frac{\varepsilon^{p}h}{p+2}\sum_{j \in \mathcal{T}_M} \left[|v^n_j|^{p+2}+ |v^{n+1}_j|^{p+2}\right] \equiv \mathcal{E}^0, \quad n \geq 0.\\
\end{split}
\end{equation}
\end{lemma}

\subsection{Main results}
Motivated by the analytical results and the assumptions on the NKGE \eqref{eq:21}, we can make the assumptions on the exact solution $v$ of the oscillatory NKGE \eqref{eq:51}: 
\begin{equation*}
(B)
\begin{split}
v \in \ &  C([0, T_0]; W_p^{6, \infty})  \cap C^2([0, T_0]; W^{4, \infty}) \cap C^4([0, T_0]; W^{2, \infty}), \\ &\quad\left\|\frac{\partial^{r+q}}{\partial s^r \partial x^q} u(x, t)\right\|_{L^{\infty}}  \lesssim \frac{1}{\varepsilon^{pr}},\quad 0 \leq r \leq 4,\quad 0 \leq r+q \leq 6,
\end{split}
\end{equation*}

Define the error function $\tilde{e}^n \in X_M (n \geq 0)$ as
\begin{equation}
\tilde{e}^n_j = v(x_j, s_n) - v^n_j, \quad j \in \mathcal{T}^0_M, \quad n \geq 0,
\label{eq:verrdef}
\end{equation}
where $v^n\in X_M$ is the numerical approximation of the oscillatory NKGE (\ref{eq:51}) obtained by the 4cFD methods. By taking $k = \varepsilon^p \tau$ in the 4cFD methods for the NKGE \eqref{eq:21}, we can directly get the error bounds of the 4cFD methods for the oscillatory NKGE \eqref{eq:51}.
\begin{theorem}
Under the assumption ($B$), there exist constants $h_0 > 0$ and $k_0 > 0$ sufficiently small and independent of $\varepsilon$, such that for any $0 < \varepsilon \leq 1$, when  $0 < h \leq h_0\varepsilon^{p/4}$, $0 < k \leq k_0 \varepsilon^{3p/2}$, the following two error estimates of the scheme (\ref{eq:I4th_HOE}) with (\ref{eq:vib}) and (\ref{eq:v1}) hold 
\begin{equation}
\|\tilde{e}^n\|_{l^2} +\|\delta^{+}_x \tilde{e}^n\|_{l^2} \lesssim \frac{h^4}{\varepsilon^p}+ \frac{k^2}{\varepsilon^{3p}},\quad \|v^n\|_{l^{\infty}} \leq 1 + M_0, \quad 0 \leq n \leq \frac{T_0}{k}.
\label{eq:eb_HOE}
\end{equation}
\label{thm:I4th_HOE}
\end{theorem}

\begin{theorem}
Under the assumption ($B$), there exist constants $h_0 > 0$ and $k_0 > 0$ sufficiently small and independent of $\varepsilon$, such that for any $0 < \varepsilon \leq 1$, when  $0 < h \leq h_0\varepsilon^{p/4}$, $0 < k \leq k_0 \varepsilon^{3p/2}$ and under the stability condition \eqref{eq:stability_HOE}, the following two error estimates of the scheme (\ref{eq:SI4th_HOE}) with (\ref{eq:vib}) and (\ref{eq:v1}) hold 
\begin{equation}
\|\tilde{e}^n\|_{l^2} +\|\delta^{+}_x \tilde{e}^n\|_{l^2} \lesssim \frac{h^4}{\varepsilon^p}+ \frac{k^2}{\varepsilon^{3p}},\quad \|v^n\|_{l^{\infty}} \leq 1 + M_0, \quad 0 \leq n \leq \frac{T_0}{k}.
\label{eq:eb_HOE}
\end{equation}
\label{thm:SI4th_HOE}
\end{theorem}

Based on the above theorems, given an accuracy bound $\delta_0 > 0$, the $\varepsilon$-scalability of the 4cFD methods for the oscillatory NKGE \eqref{eq:51} should be taken as:
\begin{equation}
h = O(\varepsilon^{p/4}\sqrt{\delta_0}) = O(\varepsilon^{p/4}),\ k= O(\varepsilon^{3p/2}\sqrt{\delta_0}) =O(\varepsilon^{3p/2}), \quad  0 < \varepsilon \leq 1.
\end{equation}
The result indicates that the 4cFD methods have better spatial resolution capacity than the FDTD methods \cite{BFY}. Also, it is useful for choosing mesh size and time step such that the numerical results are trustable.

\subsection{Numerical results of the oscillatory NKGE in the whole space}
\begin{figure}[ht!]
\centerline{\includegraphics[width = 14cm,height = 7cm]{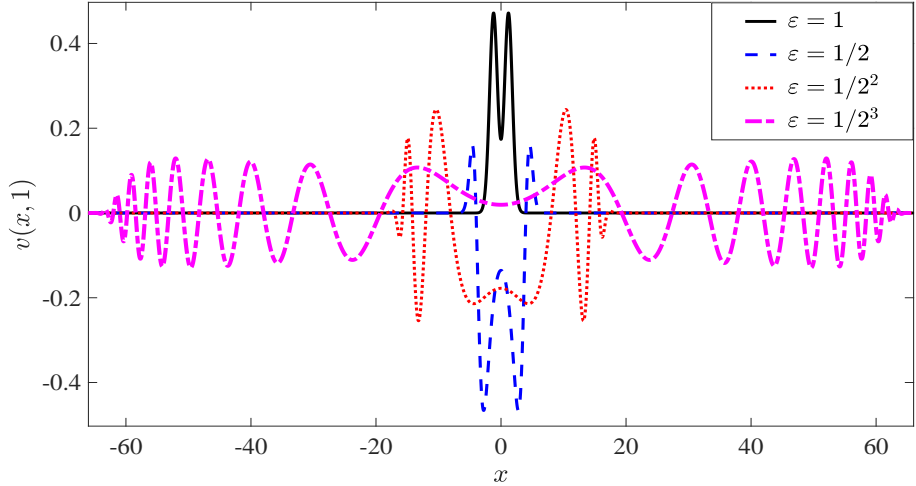}}
\caption{The solutions $v(x, 1)$ of the oscillatory NKGE \eqref{eq:550} with $d=1$, $p = 2$ and initial data \eqref{eq:ex5.1} for different $\varepsilon$.}
\label{fig:HOE_x}
\end{figure}

To avoid too much repetition, we consider the following oscillatory NKGE in $d$-dimensional ($d=1,2,3$) whole space
\begin{equation}
\begin{split}
&\varepsilon^{2p}\partial_{ss} v({\bf x,} s) - \Delta v({\bf x}, s) + v({\bf x}, s) + \varepsilon^{p} v^{p+1}({\bf x}, s) = 0,\quad {\bf{x}} \in \mathbb{R}^d,\quad s > 0,\\
&v({\bf x}, 0) = \phi({\bf x}), \quad \partial_s v({\bf x}, 0) = {\varepsilon^{-p}} \gamma({\bf x}), \quad {\bf{x}} \in \mathbb{R}^d.
\label{eq:550}
\end{split}
\end{equation}
Similar to the oscillatory NKGE \eqref{eq:50}, the solution of the oscillatory NKGE \eqref{eq:550} propagates waves with wavelength at $O(1)$ in space and $O(\varepsilon^p)$ in time, and wave speed in space at $O(\varepsilon^{-p})$. To illustrate the rapid wave propagation in space at $O(\varepsilon^{-p})$, Figure \ref{fig:HOE_x} shows the solution $v(x,1)$ of
the oscillatory NKGE \eqref{eq:550} with $d=1$, $p = 2$ and initial data
\begin{equation}
\phi(x) =  e^{-x^2} \quad \mbox{and}\quad \gamma(x) = 1/(e^{x^2} + e^{-x^2}), \qquad x\in {\mathbb R}.
\label{eq:ex5.1}	
\end{equation}

In the following, we report numerical results of the oscillatory NKGE \eqref{eq:550} with $d=1$ and $p = 1$. The initial data is chosen as \eqref{eq:ex5.1} and the bounded computational domain is taken as $\Omega_{\varepsilon} = (-4-1/\varepsilon, 4 + 1/\varepsilon)$. The `exact' solution is obtained numerically by the exponential wave integrator Fourier pseudospectral method with a very fine mesh size and a very small time step, e.g. $h_e=1/64$ and $k_e = 10^{-5}$. Denote $v^n_{h,k}$ as the numerical solution at $s=s_n$ obtained by the semi-implicit 4cFD method with mesh size $h$ and time step $k$. The errors are displayed at $s = 1$. For spatial error analysis, the time step is set as $k/\varepsilon^{2/3} = 10^{-4}$ such that the temporal error can be neglected; for temporal error analysis, we set the mesh size as $h = 1/64$ such that the spatial error can be ignored.

Tables \ref{tab:HOE_h} and \ref{tab:HOE_t} present the spatial and temporal errors for different  $0 < \varepsilon \leq 1$, respectively. From Tables \ref{tab:HOE_h} and \ref{tab:HOE_t} and additional similar numerical results not shown here for brevity, we can draw the following observations:

(i) For any fixed $\varepsilon = \varepsilon_0 > 0$, the 4cFD methods are fourth-order accurate in space and second-order accurate in time (cf. the first rows in Tables  \ref{tab:HOE_h} and \ref{tab:HOE_t}). (ii) In the highly oscillatory case, the fourth order convergence in space and second order convergence in time can be observed only when $0 < h \lesssim \varepsilon^{p/4}$ and $0 < k \lesssim \varepsilon^{3p/2}$ (cf. upper triangles above the diagonals (corresponding to $h \sim \varepsilon^{p/4}$ and $k \sim \varepsilon^{3p/2}$, and being labelled in bold letters) in Tables  \ref{tab:HOE_h} and \ref{tab:HOE_t}), which again confirm our error estimates. In summary, our numerical results confirm our rigorous error estimates and show that they are sharp.

\begin{table}[ht!]
\renewcommand{\arraystretch}{1.5}
\caption{Spatial errors of the semi-implicit 4cFD \eqref{eq:SI4th_HOE} for the NKGE \eqref{eq:550} with $d=1$, $p =1$ and initial data \eqref{eq:ex5.1}}
\centering
\begin{tabular}{ccccc}
\hline
$e_{h,k_e}(s=1)$ &$h_0 = 1/2$ & $h_0/2$ &$h_0/2^2$ & $h_0/2^3$   \\
\midrule[0.8pt]
$\varepsilon_0 = 1$  & 1.28E-2 & \bf{8.12E-4} & 5.14E-5 & 3.23E-6    \\
Order & - & \bf{3.98} & 3.98 & 3.99  \\
$\varepsilon_0 / 2^4 $  & 1.48E-1 & 1.21E-2 & \bf{8.04E-4} & 5.03E-5  \\
Order & - & 3.61 & \bf{3.91} & 4.00  \\
$\varepsilon_0 / 2^8 $  & 6.80E-1 & 1.55E-1 & 1.22E-2 & \bf{8.13E-4}  \\
Order & - & 2.13 & 3.67& \bf{3.91}  \\
\hline
\end{tabular}
\label{tab:HOE_h}
\end{table}

\begin{table}[ht!]
\renewcommand{\arraystretch}{1.5}
\caption{Temporal errors of the semi-implicit 4cFD \eqref{eq:SI4th_HOE} for the NKGE \eqref{eq:550} with $d=1$, $p =1$ and initial data \eqref{eq:ex5.1}}
\centering
\begin{tabular}{ccccccc}
\hline
$e_{h_e,k}(s=1)$ &$k_0 = 0.1 $ & $k_0/2 $ &$k_0/2^2 $ & $k_0/2^3 $ & $k_0/2^4$ & $k_0/2^5$ \\
\midrule[0.8pt]
$\varepsilon_0 = 1$ &\bf{2.25E-2} & 6.00E-3 & 1.54E-3 & 3.91E-4 & 9.83E-5 & 2.46E-5  \\
Order & \bf{-} & 1.91 & 1.96 & 1.98 & 1.99 & 2.00  \\
$\varepsilon_0 /2^{2/3}$ & 8.56E-2 & \bf{2.38E-2} & 6.18E-3 & 1.57E-3 & 3.95E-4 & 9.91E-5 \\
Order & - & \bf{1.85} & 1.95 & 1.98 & 1.99 & 1.99 \\
$\varepsilon_0 /2^{4/3}$ & 3.04E-1  & 8.92E-2 & \bf{2.35E-2} &  5.99E-3 & 1.51E-3 & 3.78E-4 \\
Order & - & 1.77 & \bf{1.92} & 1.97 & 1.99 & 2.00 \\
$\varepsilon_0 /2^{6/3}$ & 9.50E-1  & 3.31E-1 & 9.17E-2 & \bf{2.35E-2} & 5.92E-3 & 1.48E-3\\
Order & - & 1.52 & 1.85& \bf{1.96} & 1.99 & 2.00 \\
$\varepsilon_0 /2^{8/3}$ & 1.97E+0  & 1.02E+0 & 3.48E-1 & 9.33E-2 & \bf{2.36E-2} & 5.89E-3 \\
Order & - & 0.95 & 1.55 & 1.90 &\bf{1.98} & 2.00 \\
\hline
\end{tabular}
\label{tab:HOE_t}
\end{table}


\section{Conclusion}

The fourth-order compact finite difference (4cFD) methods were adapted to solve the nonlinear Klein-Gordon equation (NKGE), while the nonlinearity  strength is characterized by $\varepsilon^p$ with a constant $p \in \mathbb{N}^+$ and a dimensionless parameter  $\varepsilon \in (0, 1]$. Rigorous error bounds were established for the long time dynamics of the NKGE up to the time at $O(\varepsilon^{-p})$. The error bounds depend explicitly on the mesh size $h$, time step $\tau$ and the small parameter $\varepsilon \in (0, 1]$, which indicate the spatial/temporal resolution capacity of the 4cFD methods. Based on the error bounds, in order to obtain `correct' numerical solution of the NKGE up to the long time at $O(\varepsilon^{-p})$, the $\varepsilon-$scalability (or meshing strategy requirement) of the 4cFD methods should be taken as : $h = O(\varepsilon^{p/4})$ and $\tau = O(\varepsilon^{p/2})$, which has better spatial resolution capacity than the classical second order central difference methods. In addition, the 4cFD methods were extended to an oscillatory NKGE and the error bounds on a fixed time were obtained straightforwardly. Numerical results were reported to confirm our error bounds and demonstrate that they are sharp.

\section*{Acknowledgements}
The author would like to specially thank Professor Weizhu Bao for his valuable suggestions and comments.

%
%


\begin{thebibliography}{}
%
%

	 \bibitem{BC2}
	{W. Bao, Y. Cai},
	Optimal error estimates of finite difference methods for the Gross-Pitaevskii equation with angular momentum rotation,
	Math. Comp. vol. 82 (2013) pp. 99-128.
	
         \bibitem{BC1}
	{W. Bao, Y. Cai},
	Uniform error estimates of finite difference methods for the nonlinear Schr{\"o}dinger equation with wave operator,
	SIAM J. Numer. Anal. vol. 50 (2012) pp. 492-521.
		
         \bibitem{BCJY}
         {W. Bao, Y. Cai, X. Jia, J. Yin},
         Error estimates of numerical methods for the nonlinear Dirac equation in the nonrelativistic limit regime,
         Sci. China Math. vol. 59 (2016) pp. 1461-1494.

	\bibitem{BCZ}
	{W. Bao, Y. Cai, X. Zhao},
	A uniformly accurate multiscale time integrator pseudospectral method for the Klein-Gordon equation in the nonrelativistic limit regime,
	SIAM J. Numer. Anal. vol. 52 (2014) pp. 2488-2511.
	
	 \bibitem{BD}
	{W. Bao, X. Dong},
	Analysis and comparison of numerical methods for the Klein-Gordon equation in the nonrelativistic limit regime,
	Numer. Math. vol. 120 (2012) pp. 189-229. 
	
         \bibitem{BDZ0}
         {W. Bao, X. Dong, X. Zhao},
         An exponential wave integrator pseudospectral method for the Klein-Gordon-Zakharov system,
         SIAM J. Sci. Comput. vol. 35 (2013) pp. A2903-A2927.

         \bibitem{BDZ}
	{W. Bao, X. Dong, X. Zhao},
	Uniformly accurate multiscale time integrators for highly oscillatory second order differential equations,
	J. Math. Study vol. 47 (2014) pp. 111-150.
	
        \bibitem{BFY}
        {W. Bao, Y. Feng, W. Yi},
        Long time error analysis of finite difference time domain methods for the nonlinear Klein-Gordon equation with weak nonlinearity,
        Commun. Comput. Phys. vol. 26 (2019) pp. 1307-1334.
        
         \bibitem{BS1}
	{W. Bao, C. Su},
	Uniform error bounds of a finite difference method for the Klein-Gordon-Zakharov system in the subsonic limit regime,
	Math. Comp. vol. 87 (2018) pp. 2133-2158.
        
        	\bibitem{BS}
	{W. Bao, C. Su},
	Uniform error estimates of a finite difference method for the Klein-Gordon-Shr\"odinger system in the nonrelativistic and massless limit regimes,
	Kinet. Relat. Mod. vol. 11 (2018) pp. 1037-1062.
     	
	\bibitem{BY}
        {W. Bao, L. Yang},
        Efficient and accurate numerical methods for the Klein-Gordon-Schr\"odinger equations,
         J. Comput. Phys. vol. 225 (2007) pp. 1863-1893.
        
         \bibitem{BZ}
	{W. Bao, X. Zhao},
	A uniformly accurate (UA) multiscale time integrator Fourier pseudospectral method for the Klein-Gordon-Schr\"odinger equations in the nonrelativistic limit regime,
	Numer. Math. vol. 135 (2017) pp. 833-873.
	
	\bibitem{BZ2}
	{W. Bao, X. Zhao},
	Comparison of numerical methods for the nonlinear Klein-Gordon equation in the nonrelativistic limit regime, 
	J. Comput. Phys. vol. 398 (2019) article 108886.
	
	\bibitem{BP}
	{P. Brenner},
	On the existence of global smooth solutions of certain semi-linear hyperbolic equations,
	Math. Z. vol. 167 (1979) pp. 99-135. 
	
	\bibitem{BV}
	{P. Brenner, W. von Wahl},
	Global classical solutions of nonlinear Klein-Gordon equations,
	Math. Z. vol. 176 (1981) pp. 87-121.
	
	\bibitem{CWG}
	{Q. Chang, G. Wang, B. Guo},
	Conservative scheme for a model of nonlinear dispersive waves and its solitary waves induced by boundary motion,
	J. Comput. Phys. vol. 93 (1991) pp. 360-375.
	
	\bibitem{CCLM}
	{P. Chartier, N. Crouseilles, M. Lemou, F. M\'ehats},
	Uniformly accurate numerical schemes for highly oscillatory Klein-Gordon and nonlinear Schr\"odinger equations,
	Numer. Math. vol. 129 (2015) pp. 211-250.
	
	\bibitem{CHL}
	{D. Cohen, E. Hairer, Ch. Lubich},
	Conservation of energy, momentum and actions in numerical discretizations of non-linear wave equations,
	Numer. Math. vol. 110 (2008) pp. 113-143.
	
	\bibitem{CHL1}
	{D. Cohen, E. Hairer, Ch. Lubich},
	Long-time analysis of nonlinearly perturbed wave equations via modulated Fourier expansions,
	Arch. Ration. Mech. Anal. vol. 187 (2008) pp. 341-368.
	
	\bibitem{DG}
	{M. Dehghan, A. Ghesmati},
	Application of the dual reciprocity boundary integral equation technique to solve the nonlinear Klein–Gordon equation,
	Comput. Phys. Commun. vol. 181 (2010) pp. 1410-1418. 
	
	\bibitem{DMA}
	{M. Dehghan, A. Mohebbi, Z. Asghari},
	Fourth-order compact solution of the nonlinear Klein-Gordon equation,
	Numer. Algorithms vol. 52 (2009) pp. 523-540.
		
	\bibitem{D2}
	{J.-M. Delort},
	On long time existence for small solutions of semi-linear Klein-Gordon equations on the torus,
	J. Anal. Math. vol. 107 (2009) pp. 161-194.
		
	\bibitem{DS}
	{J.-M. Delort, J. Szeftel},
	Long-time existence for small data nonlinear Klein-Gordon equations on tori and spheres,
	Int. Math. Res. Not. vol. 37 (2004) pp. 1897-1966.
		
	\bibitem{DEGM}
	{R.K. Dodd, J.C. Eilbeck, J.D. Gibbon, H.C. Morris},
	Solitons and nonlinear wave equations,
	Academic Press, New York, 1984.
	
	\bibitem{DXZ}
	{X. Dong, Z. Xu, X. Zhao},
	On time-splitting pseudospectral discretization for nonlinear Klein-Gordon equation in nonrelativistic limit regime,
	Commun. Comput. Phys. vol. 16 (2014) pp. 440-466. 
		
	\bibitem{FZ}
	{D. Fang, Q. Zhang},
	Long-time existence for semi-linear Klein-Gordon equations on tori,
	J. Differential Equations vol. 249 (2010) pp. 151-179.
	
	\bibitem{FV}
	{H. Feshbach, F. Villars},
	Elementary relativistic wave mechanics of spin 0 and spin 1/2 particles,
	Rev. Modern Phys. vol. 30 (1958) pp. 24.
	
	\bibitem{G}
	{T.H. Gronwall},
	Note on the derivatives with respect to a parameter of the solutions of a system of differential equation,
	Ann. Math. vol. 20 (1919) pp. 292-296. 
	
	\bibitem{JV}
	{S. Jim\'enez, L. V\'azquez},
	 Analysis of four numerical schemes for a nonlinear Klein-Gordon equation,
	Appl. Math. Comput. vol. 35 (1990) pp. 61-94.
	
	
         \bibitem{KT}
	{M. Keel, T. Tao},
	Small data blow-up for semilinear Klein-Gordon equations,
	Amer. J. Math. vol. 121 (1999) pp. 629-669.
	
	 \bibitem{K}
	{S. Klainerman},
	Global existence of small amplitude solutions to nonlinear Klein-Gordon equations in four space-time dimensions,
	Comm. Pure Appl. Math. vol. 38 (1985) pp. 631-641. 

	\bibitem{L}
	{W. Liao},
	An implicit fourth-order compact finite difference scheme for one-dimensional Burgers’ equation,
	Appl. Math. Comput. vol. 206 (2008) pp. 755-764. 
	
	\bibitem{LZK}
	{W. Liao, J. Zhu, A.Q.M. Khaliq},
	A fourth-order compact algorithm for nonlinear reaction–diffusion equations with Neumann boundary conditions,
	Numer. Methods Partial Differential Equations vol. 22 (2006) pp. 600-616.
	
	\bibitem{LLG}
	{Y. Luo, X. Li, C. Guo},
	Fourth-order compact and energy conservative scheme for solving nonlinear Klein-Gordon equation,
	Numer. Methods Partial Differential Equations vol. 33 (2017) pp. 1283-1304.
	
	\bibitem{M}
	{S. Machihara},
	The nonrelativistic limit of the nonlinear Klein-Gordon equation,
	Funkcial. Ekvac. vol. 44 (2001) pp. 243-252.
	
	\bibitem{MAD}
	{A. Mohebbi, M. Abbaszadeh, M. Dehghan},
	A high-order and unconditionally stable scheme for the modified anomalous fractional sub-diffusion equation with a nonlinear source term,
	J. Comput. Phys. vol. 240 (2013) pp. 36-48. 
	
	\bibitem{MAD2}
	{A. Mohebbi, M. Abbaszadeh, M. Dehghan},
	High-order difference scheme for the solution of linear time fractional Klein--Gordon equations,
	Numer. Methods Partial Differential Equations vol. 30 (2014) pp. 1234-1253.
			
	\bibitem{SJJ}
	{J.J. Sakurai},
	Advanced Quantum Mechanics,
	Addison-Wesley, New York, 1967.
	
	\bibitem{SV}
	{W. Strauss, L. V\'azquez},
	Numerical solution of a nonlinear Klein-Gordon equation,
	J. Comput. Phys. vol. 28 (1978) pp. 271-278.
	
	\bibitem{T}
	{V. Thom{\'e}e},
	Galerkin finite element methods for parabolic problems,
	Springer, Berlin, 1997.
	
	\bibitem{VW}
	{W. von Wahl},
	Regular solutions of initial-boundary value problems for linear and nonlinear wave-equations. II,
	Math. Z. vol. 142 (1975) pp. 121-130. 
	
         \bibitem{WGX}
	{T. Wang, B. Guo, Q. Xu},
	Fourth-order compact and energy conservative difference schemes for the nonlinear Schr\"odinger equation in two dimensions,
	J. Comput. Phys. vol. 243 (2013) pp. 382-399.
              
        	\bibitem{W}
	{A.M. Wazwaz},
	New travelling wave solutions to the Boussinesq and the Klein-Gordon equations,
	Commun. Nonlinear Sci. Numer. Simulat. vol. 13 (2008) pp. 889-901. 
	
\end{thebibliography}


\end{document}